%% file: quad.tex
\documentclass{article}

\usepackage{amssymb}
\usepackage{amsmath}

\input pre.tex

\def\calQ{{\cal Q}}
\newtheorem{corollary}{Corollary}

\def\hatB{{\hat B}}
\def\outdeg{{\rm outdeg}}
\def\hxi{{\hat\xi}}
\def\tot#1{{\,\stackrel{#1}\to\,}}

\begin{document}

\title{Local structure of random quadrangulations}
\author{Maxim Krikun}
\maketitle
\begin{abstract}
This paper is an adaptation of a method used in \cite{K} to the model of random
quadrangulations. We prove local weak convergence of uniform measures on
quadrangulations 
and show that the local growth of quadrangulation is
governed by certain critical time-reversed branching process
and the rescaled profile converges to the reversed continuous-state branching process.
As an intermediate result we derieve a biparametric generating function
for certain class of quadrangulations with boundary.
\end{abstract}

\section{Introduction}
We consider the set of all finite rooted quadrangulations as a metric space
with distance between two quadrangulation defined by
\[ d(Q_1,Q_2) = \inf\B\{ \f1{1+R} \B| R: B_R(Q_1)=B_R(Q_2) \B\},
\]
where $B_R(Q)$ denotes the ball of radius $R$ around the root,
and denote the completion of this space by $\calQ$.
Elements of $\calQ$ other than finite quadrangulations are, by definition,
\defined{infinite quadrangulations}.

\begin{theorem}\label{T1}
The sequence $\mu_N$ of probability measures uniform on quadrangulations 
with $N$ faces converges weakly to a probability measure $\mu$
with support on infinite quadrangulations.
\end{theorem}
The measure $\mu$ defines certain random object -- a \defined{uniform infinite
quadrangulation}, and we are interested in local properties of this object.
We show that under distribution $\mu$ for each $R$ there exists a cycle 
$\gamma_R$, consisting of vertices at distance $R$ from the root and 
square diagonals between them, such that $\gamma_R$ separates the root 
from the infinite part of quadrangulation.
Denote by $|\gamma_R|$ the length of such cycle.
\begin{theorem}\label{T2}
$|\gamma_R|$ is a Markov chain with transition probabilities given by
\[ \P\B\{ |\gamma_{r+n}|=k \B| |\gamma_{r}|=l \B\} 
    = \f{[t^k]F(t)}{[t^l]F(t)} \cdot \P\{ \xi_n = l | \xi_0 = k \},
\]
where $\xi$ is a critical branching process with offspring generating function
\[ \phi(t) = \f1{2t} \B( \sqrt{(t-9)(t-1)^3} - 3 + 6t - t^2\B), \]
and 
\[ F(t) = \f34 \B( \sqrt{\f{9-t}{1-t}} - 3 \B) \] 
is the generating function of it's stationary measure.
\end{theorem}

\begin{corollary}\label{C1}
The random variable $2 |\gamma_R|/R^2$ converges in distribution to the $\Gamma(3/2)$ law.
\end{corollary}
Knowing that the properly rescaled branching processes converge to 
continuous-time branching processes, it is natural to look for 
the continuous-time limit of the rescaled profile $|\gamma_{tR}|/R^2$.
An exact statement is provided by Theorem~\ref{T3} in Section~\ref{sec.cbp}.


Finally, we propose the following conjecture.
\begin{conjecture}\label{H1}
Let $\ell(R)$ be the cycle of minimal length in a uniform infinite
quadrangulation $Q$ that separates the ball $B_R(Q)$ from the
infinite part of quadrangulation.
Then it's length $|\ell(R)|$ is linear in $R$ as $R\to\infty$ 
\end{conjecture}
In \ref{sec.l} we show how the upper bound follows from the branching
construction used in proof of Theorem~\ref{T2}.
The lower bound remains an open question.

\section{Some facts on quadrangulations}\label{sec.d}

\subsection{Definitions}
Consider a finite planar graph embedded into the sphere, such that each 
component of the complement to the graph is homeomorphic to a disk.
A \defined{planar map} is an equivalence class of such embedded graphs 
with respect to orientation-preserving homeomorphisms of the sphere.

A planar map is \defined{rooted} if a directed edge, called the 
\defined{root}, is specified. 
A rooted planar map has no nontrivial automorphisms.
We will refer to the tail vertex of the root as \defined{root vertex},
and to distance from any vertex of a map to this vertex as 
\defined{distance to the root}.

\defined{Quadrangulation} is a rooted planar map such that all it's faces
are squares, i.e. it's dual graph is four-valent.
Note that every quadrangulation is bipartite (this follows from the
fact that any subset of faces in a quadrangulation is necessarily
bounded by an even number of edges).

In the following we will distinguish two type of faces based on the 
distances from the vertices around the face to the root: 
these distances are either $(R-1,R,R+1,R)$ or $(R-1,R,R-1,R)$ for 
some $R$. Since every quadrangulation is bipartite, there are no $(R,R)$
edges.

\defined{Quadrangulation with a boundary} is rather self-evident notion;
formally this is a planar map with all faces being squares except one 
distinguished face which can be an arbitrary even-sided polygon.
When drawing the quadrangulation it is convenient to represent this 
distinguished face by the infinite face. This face is then excluded from 
"faces" of quadrangulation and is referred to as "boundary".

We say that a quadrangulation has \defined{simple boundary}, if
all vertices of the boundary are distinct (i.e. no vertex is met twice
when walking around the boundary), and every second vertex has degree two
(see \figref{simple}).

\putfigure{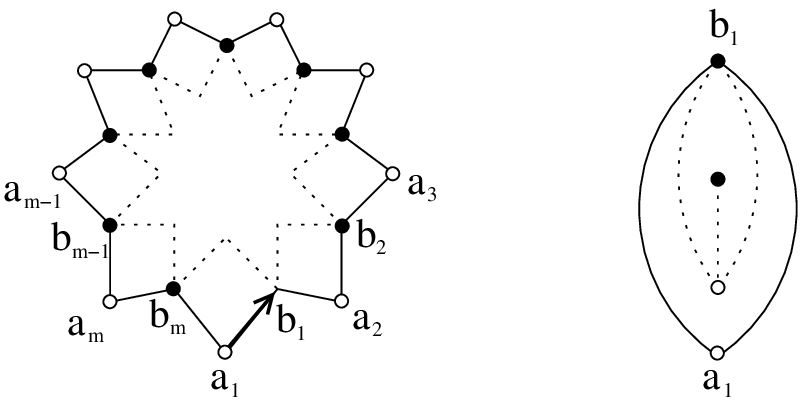}{Typical quadrangulation with simple boundary
           and an example with minimal boundary length.}

\subsection{Some enumeration results}\label{sec.enum}

Let $C(N)$ be the number of rooted quadrangulations with $N$ faces, and
let $C(N,m)$ be the number of rooted quadrangulations with $N$ faces and
with simple boundary of length $2m$, such that the root is located 
on the boundary and root vertex has degree two.

We will need the following enumeration results
(see section \ref{sec.p} for details)
\begin{equation}\label{eq.qx}
 q(x) = \sum_{N=0}^\infty C(N) x^N 
        = \f43 \f{2\sqrt{1-12x}+1}{(\sqrt{1-12x}+1)^2},
\end{equation} 
\begin{eqnarray}  
U(x,y) 
 &=& \sum_{N,m} C(N,m) x^N y^m \nonumber\\
 &=& \f12\B(y-xy^2-1 + \sqrt{ y^2-2xy^3-2y + 4xy q(x) + (xy^2-1)^2}\B)
\label{eq.Uxy}
\end{eqnarray}
The function $U(x,y)$ is analytic around $(0,0)$ and it's first singularity 
in $x$ for small $y$ coincides with the singularity of $q(x)$, 
i.e. $x_0 = 1/12$.
From the expansion near this point
\begin{equation}\label{eq.uexp} 
  U(x,y)\B|_{x=x_0-t} = A(y) + A_1(y) t + B(y) t^{3/2} + O(t^2), 
\end{equation}
where
\[ A(y) = \f1{24}\sqrt{(y-18)(y-2)^3} - \f12 + \f{y}2 - \f{y^2}{24}, \]
\[ A_1(y) = \f{y^2}2 + \f{y}2 \f{y^2-10y-32}{\sqrt{(18-y)(2-y)}}, \qquad
   B(y) = \f{64 \sqrt3\cdot y }{ \sqrt{(y-18)(y-2)^3} }, \]
one finds the asymptotic of $C(N,m)$ as $N\to\infty$:
\begin{equation}\label{eq.asympt} 
  C(N,m) = \f{b(m)}{\Gamma(\f32)} N^{-5/2} x_0^{-N} \B( 1+O(N^{-1/2})\B), 
   \qquad
   b(m) = [y^m] B(y).
\end{equation}
Note also that $[y]U(x,y)=x q(x)$, thus
\begin{equation}\label{eq.asympt2}
 C(N) = C(N+1,1) 
      = \f{b(1)}{x_0 \Gamma(\f32)} N^{-5/2} x_0^{-N} \B( 1+O(N^{-1/2})\B). 
\end{equation}

\subsection{Basic probabilities}

First let us specify more exactly the definition of ball $B_R(Q)$. 
Given a rooted quadrangulation $Q$, $B_R(Q)$ consists of all faces that have
at least one vertex at distance strictly less than $R$ from the root.
With this definition there are only faces of type $(R-1,R,R+1,R)$ at
the boundary of $B_R(Q)$.

Say we want to compute the probability for a uniformly distributed $N$-faced 
quadrangulation $S_N$ to have a particular root neighbourhood $K$. Suppose that
$K$ has $n$ faces and, for simplicity, a connected boundary of length $2m$,
so that $K$ is a quadrangulation with simple boundary.
Take any other quadrangulation $S$ with simple boundary of the same length.
We can glue $K$ and $S$ as follows: 
\begin{itemize}
\item{} cut $m$ half-squares around the boundary of $K$, so that 
        the resulting map $K'$ is bounded by $m$ diagonals;
\item{} repeat the same for $S$, obtaining a map $S'$ bounded by $m$ diagonals;
\item{} identify the boundaries%
\footnote{ There are $m$ different ways to do this, but since $K$ is rooted
           we can choose (in some deterministic way) one of diagonals on the
           boundary of $K'$ and require that when glueing it is identified 
           with the diagonal of the rooted face of $S$; 
           this way the ambiguity is eliminated. }
of $K'$ and $S'$ (see \figref{glue}).
\end{itemize}
The resulting map $(K'+S')$ has $m$ faces less than $K$ and $S$ had together, 
i.e. if $S$ has $N-n+m$ faces, $(K'+S')$ would have $N$ faces.

\putfigure{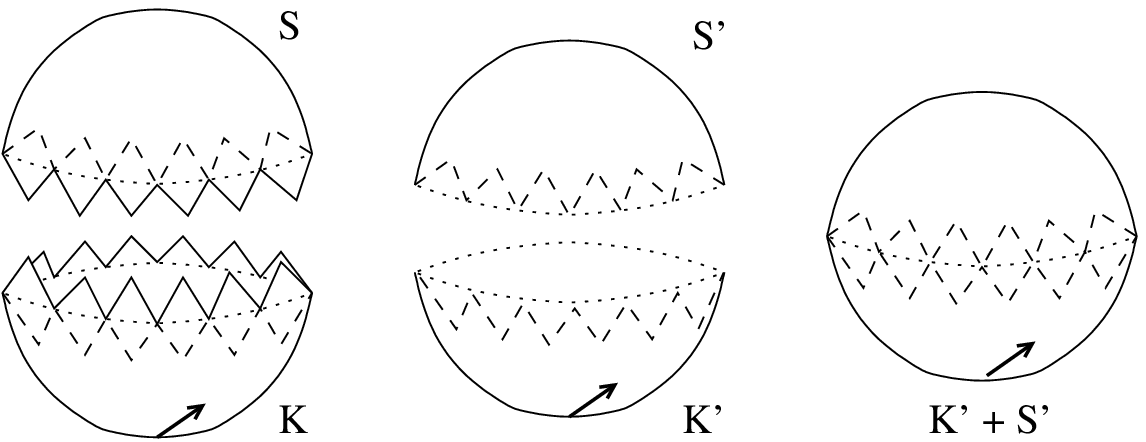}{Glueing two quadrangulations with simple boundary.}

It's easy to see that the process on \figref{glue} is reversible.
Indeed, take a quadrangulation $Q$ with root neighbourhood $K$, cut it in two 
along the boundary of $K'$, and add $m$ half-squares to each part. 
This will give $K$ and some quadrangulation with simple boundary, which
then can be used to reconstruct $Q$.

Thus for each of $C(N-n+m, m)$ maps $S$ with $N-n+m$ faces we get a different
$N$-faced quadrangulation, and every $N$-faced quadrangulation with root 
neighbourhood $K$ is obtained this way.
In other words
\[ \P\{ B_R(S_N)=K \} = \f{C(N-n+m, m)}{C(N)}. 
\]
Combining this with asymptotics \eqref{asympt}, \eqref{asympt2} we get
\begin{lemma}\label{L1}
Given a quadrangulation $K$ with $n$ faces and simple boundary of length $2m$,
such that $B_R(Q)=K$ for some $Q$,
\begin{equation}\label{eq.psn1}
 \lim_{N\to\infty} \P\{ B_R(S_N) = K \} = \f{ b(m) }{ b(1) } x_0^{n-m+1},
\end{equation}
where $S_N$ denotes uniformly distributed random quadrangulation with $N$ faces.
\end{lemma}

\putfigure{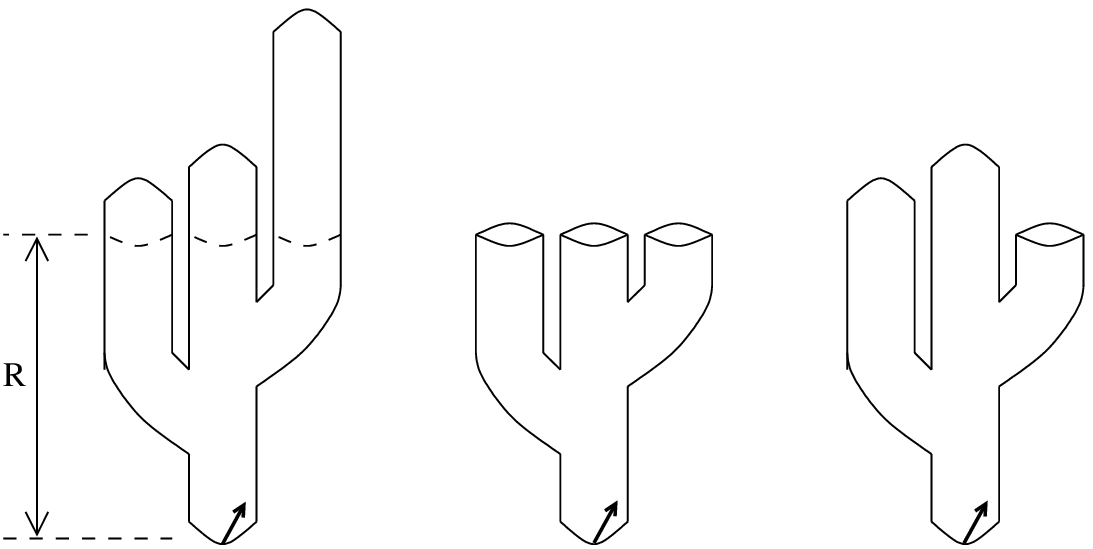}{Quadrangulation $Q$, ball $B_R(Q)$ and hull $\hatB_R(Q)$.}

In general case, however, the boundary of $B_R(Q)$ may have multiple disjoint
components (\figref{cut}, middle).
Following the same reasoning as above and assuming that $K$ has $k$ 
boundary components ("holes") of length $2 m_1,\ldots, 2 m_k$ and $n$ faces,
we'll get the following formula
\begin{equation}\label{eq.pksn}
\P\{ K = B_R(S_N) \} = \f{1}{C(N)} \sum_{N_1+\ldots+N_k = N -n}
                       \,\,\prod_{j=1}^k C(N_j+m_j, m_j).
\end{equation}
Here we count all possible ways to "fill" the $k$ holes in $K$ using 
quadrangulations with appropriate boundary length; $N_j$ is the number
of internal faces in quadrangulation used to fill $j$th hole.

Due to the factor $N^{-5/2}$ in asymptotics \eqref{asympt},
for large $N$ the only significant terms in sum~\eqref{pksn} 
are those where one of $N_j$ has order $N$, while all others are finite.
This means that in a large random quadrangulation $S_N$ conditioned to 
$B_R(S_N)=K$, with high probability only one of the "holes" in $K$ contains 
the major part of the quadrangulation (we could calculate exact probabilities
here, but this is not necessary).

Such observation motivates the following definition:
given quadrangulation $Q$, take the ball $B_R(Q)$ and glue all 
but the largest components of the complement $Q\backslash B_R(Q)$ 
back to the ball.%
\footnote{ if there are multiple components with maximal size, let us chose
           one of them in some deterministic way; details are not important
		   to us since such situations have small probability for large $N$. }
The resulting map is called the \defined{$R$-hull} of quadrangulation $Q$,
and is denoted by $\hatB_R(Q)$.

Unlike the boundary of the ball, the boundary of $\hatB_R(Q)$ is always
connected (see \figref{cut}, right), but similarly to $B_R(Q)$ there 
there are only faces of type $(R,R+1,R,R-1)$ at the boundary of $\hatB_R(Q)$,
thus the hull is a quadrangulation with simple boundary.

Limiting probabilities for the hull and exactly the same as for the ball:
\begin{lemma}\label{L2}
Given a quadrangulation $K$ with $n$ faces and boundary of length~$2m$,
such that $\hatB_R(Q)=K$ for some $Q$,
\begin{equation}\label{eq.psn2}
 \lim_{N\to\infty} \P\{ \hatB_R(S_N) = K \} = \f{ b(m) }{ b(1) } x_0^{n-m+1}
\end{equation}
where $S_N$ denotes uniformly distributed random quadrangulation with $N$ faces.
\end{lemma}
\proof
The proof is essentially the same as for the ball with single hole.
Given $K$, every quadrangulation $(K'+S')$, obtained by glueing $K$ and any
quadrangulation $S$ with simple boundary of length $2m$, has the same $R$-hull
$\hatB_R(K'+S')=K$ as soon as the number of faces in $S$ is large enough
(say larger than $n$).
Thus for $N>2n$
\[ \P\{ \hatB_R(S_N) = K \} = \f{C(N-n-m,m)}{C(N)} \]
and the limit \eqref{psn2} follows. 
\eop

\subsection{A note on convergence of measures}
The limiting probabilities~\eqref{psn2} define a measure $\mu$ on $\calQ$,
such that for all~$K$ and~$R$
\[ \mu_N( Q :\, B_R(Q)=K) \to \mu(Q:\, B_R(Q)=K) 
   \quad\mbox{as $N\to\infty$}.
\]
However, since $\calQ$ is not compact, the existence of this limit
does not, by itself, imply weak convergence of $\mu_N$ to $\mu$.
For the weak convergence to follow one has to show that $\mu$ is indeed a
probability measure. See \cite[section 1.2]{AS} for detailed discussion of this question.

In the next section we will evaluate the sum of limiting
probabilities \eqref{psn2} over all possible $R$-hulls $K$ and show that for
each $R$ this sum equals one. This will prove Theorem~\ref{T1}.

\section{Quadrangulation and branching process}\label{sec.q}

\subsection{Hull decomposition}
Consider $K$ such that $K=\hatB_R(Q)$ for some quadrangulation $Q$.
If $Q$ is large enough (e.g. if the number of faces in $Q$ is at least
twice that of $K$) then
\[ \hatB_R(Q) \supseteq \hatB_{R-1}(Q) \supseteq \ldots \supseteq \hatB_1(Q), \]
and this sequence doesn't actually depend on $Q$.

As noticed earlier, the hull has simple boundary. 
Let us denote the vertices of the boundary of $\hatB_R(Q)$ by 
$(a_1, b_1, \ldots, a_m, b_m)$, as on \figref{simple}, starting 
from some arbitrarily chosen vertex and so that all $a_i$'s are situated 
at distance $R+1$ from the root, and all $b_i$'s at distance $R$ from the root.

Let $\gamma_R$ be the cycle consisting of vertices $b_1,\ldots,b_m$ and
square diagonals between them. Define cycles $\gamma_{R-1}, \ldots, \gamma_1$
similarly.
A \defined{layer} $L_R$ is a part of quadrangulation contained between 
cycles $\gamma_{R-1}$ and $\gamma_{R}$. 
It turns out that the layer has very simple structure:
\begin{itemize}
\item{} each edge $(b_{i-1},b_i)$ of it's \defined{upper boundary} $\gamma_R$ 
        is a diagonal in a square that touches it's \defined{lower boundary}
		$\gamma_{R-1}$ at some point $c_i$;
\item{} points $c_i$ are cyclically ordered around $\gamma_{R-1}$ and there
        are $l_i\ge0$ edges of $\gamma_{R-1}$ between $c_{i}$ and $c_{i+1}$
		($l_i=0$ if $c_i$ and $c_{i+1}$ both refer to the same vertex).
\end{itemize}
Let us call the area $(c_i,b_i, b_{i+1}, c_{i+1})$ a \defined{block}.
A layer is uniquely (up to rotation related to choice of vertex $a_1$) 
characterized by a sequence of blocks.
\putfigure{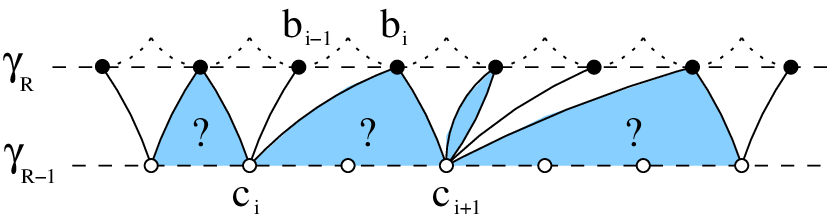}{Layer structure. Contents of filled areas is unknown.}

The internal structure of a block is not as simple -- it can contain 
arbitrary large subquadrangulation, which can have vertices at distance 
more than $R$ from the root. This is here where the "reattached" components of
$\hatB_R(Q)\backslash B_R(Q)$ go.
Note that even with $l_i=0$ the block contents can be non-trivial 
(\figref{blocks}, right).
\putfigure{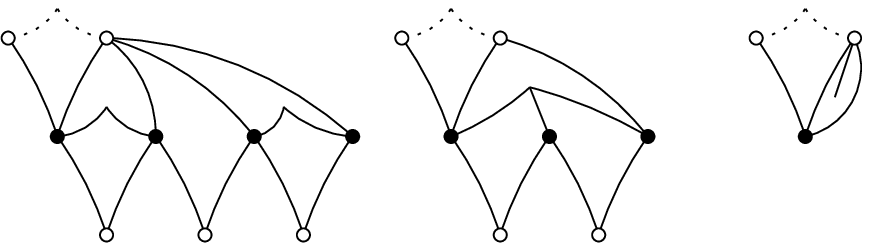}{Possible internal structure of a block.}

Fortunately there is a bijection between blocks and a class of quadrangulations 
with simple boundary counted by $C(N,m)$ in section~\ref{sec.enum}.
The block is converted to quadrangulation via the procedure illustrated 
on~\figref{convert}. Clearly, this procedure is reversible:
one has to choose the topmost vertex on the right-hand
side of \figref{convert} as the root vertex of the quadrangulation;
then the block is recovered by cutting the quadrangulation along the edge,
opposite to the root in the rooted square.
\putfigure{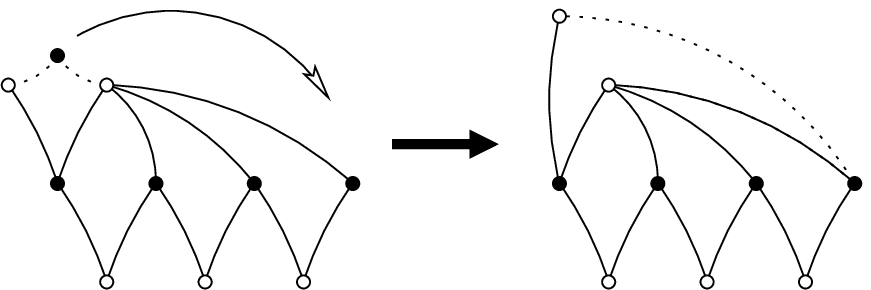}{Converting a block to a quadrangulation with simple boundary.}

To conclude: the $R$-hull consists of $R$ layers, each layer consists of
one or more blocks, and each block is essentially a quadrangulation with
simple boundary.

\subsection{Tree structure}
The layer/block representation suggests the following tree structure:
let the edges of $\gamma_r$, $r=1,\ldots, R$ be the nodes of a tree, and connect
each edge of $\gamma_r$ to the edges of $\gamma_{r-1}$ that belong to the 
same block. 
\putfigure{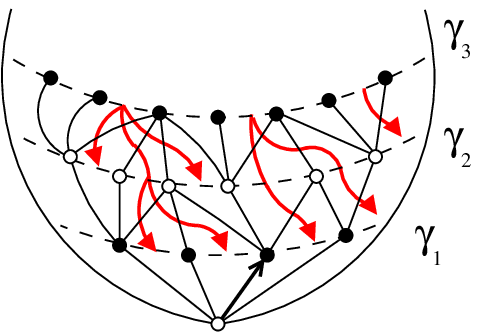}{Fragment of tree structure on quadrangulation.}

The whole hull $K$ then can be represented by a planar forest $T_K$
of height%
\footnote{note that this forest is "reversed" with respect to $K$: it starts at
          $\gamma_R$ and grows down to the root of $K$. In the following
		  we will keep using such reversed notation and will refer to nodes
		  corresponding to $\gamma_r$ as the $r$-th level of $T_K$.
		  }%
~$R$, where with each vertex one associates a quadrangulation with 
simple boundary, so that for a vertex $v_i$ of outdegree $l_i$ the associated 
quadrangulation has boundary length $2(l_i+1)$. 
Unfortunately this representation is not unique: given $T_K$ 
and it's associated quadrangulations, we can reconstruct $K$ and the 
root vertex of $K$, but not the root edge.
In order to include the information on root edge position into the tree
structure, we apply to $K$ the following modification (see \figref{extend}):
\begin{itemize}
\item{} cut $K$ along the root edge, obtaining a hole of length two;
\item{} attach a new square to the boundary of this hole;
\item{} identify two remaining edges of this square and make the 
        resulting edge a new root.
\end{itemize}
\putfigure{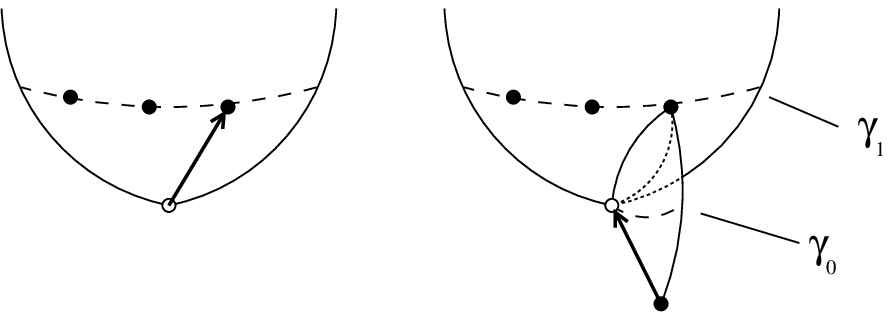}{Adding an extra square at the root.}
One diagonal of a new square has it's ends identified;
this gives an extra cycle $\gamma_0$, which always has length one.
In terms of tree structure this means that we add 
one child to some $\gamma_1$-vertex of $T_K$; 
this new vertex has no associated quadrangulation.
Call this extended forest a \defined{skeleton} of hull $K$.

Note that since there is no natural "first" edge in $\gamma_R$, the tree
structure implies only cyclic order on the trees of $T_K$.
However for convenience we will consider $T_K$ as linearly ordered,
and will keep in mind that the same tree structure can be represented
by several forests, which differ by cyclic permutation of trees.

Apart from this ambiguity the hull $K$ is uniquely characterized by 
it's skeleton and associated quadrangulations.

\subsection{Analogy with branching process}
\defined{Branching process} is a random process with discrete time.
It starts with one or more \defined{particles}, and at each step every particle
independently of the others is replaced by zero, one, or more 
\defined{child particles} according to the \defined{offspring distribution} 
$\{p_i\}_{i=0}^\infty$, which remains fixed throughout the whole process.

It is convenient to represent the trajectory of a branching process by a planar
tree (or forest, if starting from multiple particles). The probability to see 
certain trajectory tree $T$ is then a product over all vertices of 
probability for a particle to have an offspring of size equal to the 
outdegree of the vertex:
\begin{equation}\label{eq.bp} 
\P\{ T \} = \prod_{v\in T} \left. p_i \right|_{i=\outdeg(v)} 
\end{equation}
We will attempt to do the reverse: apply the theory of branching processes to
the analysis of the tree structure described above.

Say we want to compute the probability for an $R$-hull of uniformly distributed
quadrangulation $S_N$ to have a particular skeleton $T_K$.
As explained above, every such $R$-hull is obtained from $T_K$ by choosing an
appropriate set of associated quadrangulations.
On the other hand, taking for every vertex $v_i \in T_K$ with outdegree $l_i$
a quadrangulation with simple boundary of length $2(l_i+1)$ will give a 
valid $R$-hull, which has the required skeleton.
A simple calculation shows that if $i$'th associated quadrangulation has $n_i$
faces, the hull will have  $m - 1 + \sum_i (n_i-1)$ faces, 
where $m$ is half the length of hull boundary (or equivalently the number of
trees in the skeleton).

{\samepage
Combining this with Lemma~\ref{L2}, we'll get the following formula
\begin{eqnarray}
\lim_{N\to\infty} \lefteqn{ \P\{ skel_R(S_N) = T_K \} }&& \nonumber\\
  &=& \sum_{n_1,n_2,\ldots=0}^\infty
    \B(\prod_{v_i \in T_K} C(n_i, l_i+1)\B) 
	\f{b(m)}{b(1)} x_0^{\Sigma_i (n_i-1)}\nonumber\\
  &=& \f{b(m)}{b(1)} 
    \prod_{v_i \in T_K} \sum_{n_i=0}^\infty C(n_i, l_i+1) x_0^{n_i-1}\nonumber\\
  &=& \f{b(m)}{b(1)}
    \prod_{v_i \in T_K} \f1{x_0} [y^{l_i+1}]A(y),
\label{eq.pskel}
\end{eqnarray}
where $A(y)$ is the first coefficient of expansion~\eqref{uexp}.
}

The last product in~\eqref{pskel} looks similar to the product \eqref{bp}.
There is however one important difference -- product terms of \eqref{pskel} do
not define a probability distribution.

In order to make the analogy with branching process complete,
we apply the following normalization procedure: 
for each square crossed by
one of the cycles $\gamma_0, \ldots, \gamma_R$ write $y_0$ on the upper
half-square and $y_0^{-1}$ on the lower half-square.
Each block then gets an extra term $y_0^{l_i-1}$; plus there is $y_0^{m}$ on
the upper boundary of the hull, due to $m$ half-squares above $\gamma_R$,
and $y_0^{-1}$ due to one half-square below $\gamma_0$.
So we can continue \eqref{pskel} with
\begin{eqnarray}
  \ldots 
  &=& \f{b(m)}{b(1)} y_0^{m-1}
    \prod_{v_i \in T_K} \f1{x_0} y_0^{l_i-1} [y^{l_i+1}]A(y)\nonumber\\
  &=& \f{b(m)}{b(1)} y_0^{m-1} 
     \prod_{v_i \in T_K} [t^l_i] \phi(t),
\label{eq.pskel2}
\end{eqnarray}
where
\[ \phi(t) = \f{1}{t x_0 y_0^2} A(y_0 t). \]
For the Taylor coefficients of $\phi(t)$ to define a probability distribution
it has to satisfy the equation $\phi(1)=1$, which is equivalent to
$A(y_0) = x_0 y_0$. Solving this last equation we find $y_0=2$ and
\[ \phi(t) = \f1{2t} \B( \sqrt{(t-9)(t-1)^3} - 3 + 6t - t^2\B). \]

\subsection{Remaining proofs}
\proofof{Theorem~\ref{T1}}
Using \eqref{pskel}, \eqref{pskel2} let us compute the probability 
$|\gamma_R| = m$ with respect to $\mu$.
\begin{equation}\label{eq.pgm}
\P_\mu\{ |\gamma_R| = m \} 
   = \f1m \f{b(m)}{b(1)} y_0^{m-1} 
		  \sum_{T} \prod_{v_i \in T} [t^l_i] \phi(t),
\end{equation}
where the sum is taken over all forests $T$ of height $R+1$ that have
$m$ vertices on level $R$ and exactly one vertex at level $0$ (in "reversed"
notation).
The term $1/m$ appears in \eqref{pgm} because for each $R$-hull with
$|\gamma_R|=m$ there are exactly $m$ linearly ordered forests describing the
tree structure of this hull.

Let $\xi$ be a branching process with offspring generating function $\phi(t)$.
The sum in \eqref{pgm} can be interpreted as the probability for $\xi$ starting
from state $m$ at time $0$, to reach state $1$ at time $R$.
Let 
\[ F(t) = \sum_{m=0}^\infty \f{b(m)/m}{b(1)} y_0^{m-1} t^m.
\]
Then \eqref{pgm} can be rewritten as
\begin{equation}\label{eq.pgm2}
 \P_\mu\{ |\gamma_R| = m \}
   = [t^m] F(t)\cdot \P\{ \xi_{R+1}=1 | \xi_{0} = m \}.
\end{equation}	 
The $R$-step transition probabilities of a branching process are expressed via
it's offspring generating function as 
\[ \P\{ \xi_R = n | \xi_{0} = m \} = [t^n] \phi_R^m(t), \]
where $\phi_R(t)$ stands for the $R$'th iteration of $\phi(t)$.
Thus
\[
  \sum_{m=1}^\infty \P_\mu\{ |\gamma_R| = m \}
   = [t] F(\phi_R(t)). 
\]
Since $b(m)=[y]B(y)$, where $B(y)$ is a coefficient in \eqref{uexp},
we find that 
\[ F(t) = \f1{y_0 b(1)} \intl_0^{t y_0} \f{B(y)}{y}\, dy
		= \f34 \B( \sqrt{\f{9-t}{1-t}} - 3 \B)
\]
and a direct computation shows that $F(t)$ satisfies the Abel equation
\begin{equation}\label{eq.abel}
 F(\phi(t)) - F(\phi(0)) = F(t).
\end{equation}
In particular this means that $[t] F(\phi_r(t)) = [t] F(t)$ for all $r$,
and since $[t]F(t)=1$ 
\[ \sum_{m=1}^\infty \P_\mu\{ |\gamma_R| = m \} = 1. \]
But the last sum is also the sum of limiting probabilities \eqref{psn2}
over all possible $R$-hulls. This completes the proof of Theorem~\ref{T1}.
\eop

\proofof{Theorem~\ref{T2}}
To prove theorem~\ref{T2} first note that 
\[
\P_\mu \{ |\gamma_r|=l, |\gamma_{r+n}|=k \}
   = [t^k]F(t) \cdot \P\{ \xi_n = l, \xi_{r+n} = 1 | \xi_0 = k \}.
\]   
This formula is obtained by taking sum of probabilities \eqref{pskel2} 
over all skeletons with $k$ vertices at level $(r+n)$ and $l$ vertices 
at level $r$.
Now since $\xi$ is Markovian
\[ \P\{ \xi_n = l, \xi_{r+n} = 1 | \xi_0 = k \}
   = \P\{ \xi_n = l|\xi_0 = k \} \P\{ \xi_{r+n}=1 | \xi_{n} = l \},
\]
and combining this with \eqref{pgm2} we get
\begin{equation}\label{eq.rbp}
 \P\B\{ |\gamma_{r+n}|=k \B| |\gamma_{r}|=l \B\} 
    = \f{[t^k]F(t)}{[t^l]F(t)} \cdot \P\{ \xi_n = l | \xi_0 = k \}.
\end{equation}
The Abel equation~\eqref{abel} means that $F(t)$ is a generating function of 
stationary measure for process $\xi$ (see \cite{HS}, section 1.4), 
so the right-hand side of~\eqref{rbp} is indeed the transition probability for a
reversed branching process.
\eop

\proofof{Corollary~\ref{C1}}
The $R$th iteration of $\phi(t)$ is
\begin{equation}\label{eq.phiiter}
 \phi_R(t) = 1 - \f8{ \B(\sqrt{\f{9-t}{1-t}}+2R \B)^2 - 1}, 
\end{equation}
this can be verified by induction.
The distribution of $|\gamma_R|$ is given by
\[ \P\{ |\gamma_R| = m \} = [t^m] F(t) \cdot [t] \phi_R^m(t). \]
Calculating explicitly
\[ [t] \phi_R^m(t) = \f43 m (2R+3) \f{ (R^2+3R)^{m-1} }{ (R^2+3R+2)^{m+1} }, \]
\[ [t^m] F(t) = \f{3}{\sqrt{2\pi}} m^{-1/2} + O(1), \]
and putting $m = x R^2$ we find for large $R$
\[ \P\{ |\gamma_R| = xR^2 \}
   = \f8{\sqrt{2\pi}} \f{ x^{1/2} e^{-2x} }{ R^2 } + O(\f1{R^3}).
\]
This implies convergence of $2 |\gamma_R|/R^2$ to $\Gamma(3/2)$ law.

\subsection{Linear cycle}\label{sec.l}
%
%

%

The cycle $\gamma_R$ is a natural analog of circle in Euclidean geometry: 
this is a closed curve, and it's points are situated {\em exactly} at distance
$R$ from the center (root).  The relation $|\gamma_R| = O(R^2)$ is, however,
quite different from the usual $L=2\pi R$.

A natural question to ask is what happens, if we weaken restrictions, for
example, by allowing the separating cycle to contain any points at distance 
{\em at least} $R$ from the root?

It turns out that there exists such cycle with length {\em linear} in $R$.
This cycle is built as follows:
\begin{itemize}
\item{} consider all vertices of $\gamma_R$, and group together the vertices
        that have common ancestor in $\gamma_{2R}$;
\item{} in each group there is a leftmost element and a rightmost element. 
        Take a path from the leftmost up to the common ancestor, then down to
		the rightmost (such path has approximately $2R$ steps);
\item{} join these paths together to form a separating cycle $\ell_R$.
		$\gamma_{2R}$ that have non-empty offspring at $\gamma_R$.
\end{itemize}		
The length of $\ell_R$ is $2R\cdot\theta_R$, where $\theta_R$ is the number of
vertices at $\gamma_{R2}$ that have nonempty offspring at $\gamma_R$.
It remains to show that $\theta_R$ has finite distribution.
\putfigure{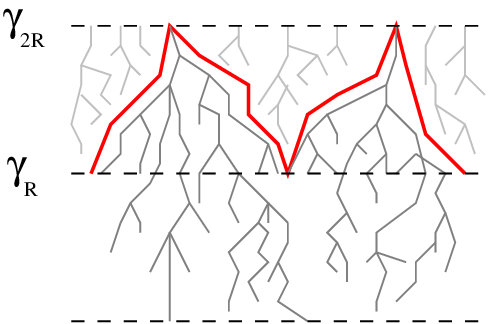}{fragment of linear cycle}

\par\noindent
From~\eqref{rbp} we have
\[ \P\{ \theta_R =k | \gamma_R = m \}
   = \sum_{n=1}^\infty \f{[t^n]F(t)}{[t^m]F(t)}
     { n \choose k } \phi_R^{n-k}(0) \cdot [t^m] (\phi_R(t)-\phi_R(0))^k,
\]
\[ \P\{ \gamma_R = m \} = [t^m] F(t) \cdot [t] \phi_R^m(t), \]
thus
\[ \E( y^{\theta_R} | \gamma_R = m )
   = \f{ [t^m] F\B( \phi_R(0) + y (\phi_R(t)-\phi_R(0)) \B) }
       { [t^m] F(t) },
\]   
\[ \E y^{\theta_R} = [t] F\B( \phi_R(0) + y (\phi_{2R}(t) - \phi_R(0)) \B).
\]
The last expression being a generating function of $\theta_R$,
a direct calculation of it's derivatives at $y=1$ shows that as $R\to\infty$ 
the first and second moments of $\theta_R$ converge to
$11/2$ and $171/2$ respectively.
This provides an upper bound for Conjecture~\ref{H1}.

\section{Convergence to continuous process}\label{sec.cbp}

\subsection{Identifying the limit}

Given the results of previous section, it is natural to expect the
convergence of rescaled process
\begin{equation}\label{eq.gatr}
\f{|\gamma_{tR}|}{R^2}, \qquad t\in[0,1] 
\end{equation}
as $R\to\infty$ in sence of weak convergence in $D[0,1]$
(the space of functions without discontinues of second kind \cite{GS}).
We will start from identifying the limit $\zeta_t$ of the rescaled 
branching process
\begin{equation}\label{eq.xitr}
\f{\xi_{tR}}{R^2}, \qquad t\in[0,1].
\end{equation}

As shown by Lamperti \cite{La}, 
any possible limit of a sequence of rescaled branching processes is a 
\defined{countinous-state branching process}, i.e. a Markov process 
on $[0,\infty)$ with right-continuous paths, 
whose transition probabilities satisfy the \defined{branching property}
\[ P_t(x+y, \cdot) = P_t(x,\cdot) * P_t(y,\cdot) \]
for all $t, x, y \ge 0$.
Every continuous-state branching process is fully characterized by it's 
Laplace transform
\[
\int_0^\infty e^{-\lambda y} P_t(x,dy) = \exp( -x  u_t(\lambda)), 
\quad \lambda\ge0,
\]
and in particular by the function $ u_t$.
To derieve the function corresponding to $\zeta_t$
we use a theorem due to Zolotarev (\cite{Z}, Theorem~7), 
\begin{theorem}[Zolotarev]\label{TZ}
Let $X_t$ be a continuous-time critical branching process with offspring 
generating function $f(t)$. If $f(1-x)$ is properly variable at zero
with index $1<1+\alpha\le 2$ in sence of Karamata, 
then the Laplace transform corresponing to the distribution 
\[ \B( \f{X_t}{\P\{X_t>0\}} \B| X_t > 0 \B) \]
converges as $t\to\infty$ to
\[ 1 - s  (1+s^\alpha)^{-1/\alpha}. \]
\end{theorem}
Thanks to the explicit formula \eqref{phiiter} we can consider a family of
functions $\{\phi_R(\cdot), R\in \mathbb{R}_+\}$ as a semigroup with generator
\[ f(t) = \lim_{R\to 0} \f{\phi_R(t) - t}{R}
        = \f12 (9-t)^{1/2} (1-t)^{3/2}.
\]
Let $\hxi$ be a continuous-time branching process with offspring generating 
function $f(t)$ (that is a particle of $\hxi$ branches into $k$ child
particles at rate $[t^k] f(t)$ in continous time).
The generating function of $\hxi_t$,
\[ F(t,x) = \E x^{\hxi_t}, \]
has to satisfy the differential equation
\[ \f\d{\d t} F = f(F), \quad F(0,x)=x \]
which is also a semigroup equation for functions $\phi_R$, i.e.
\[ F(t,x) = \phi_t(x). \]
Thus the process $\hxi_t$, taken at integer times, 
and the process $\xi_t$ are identically distributed.
Applying Theorem~\ref{TZ} to $\xi_R$ we find
\[ \lim_{R\to\infty}
   \E \B( e^{ - 2s \xi_R / R^2} \B| \xi_0=1, \xi_R>0 \B)
   = 1 - \f{s}{(1+\sqrt{s})^2}
\]
(note that the scaling factor here is $R^2/2$ rather than $R^2$).

Now note that for the branching process $\xi$ started from $\xi_0=xR^2/2$
the number of initial particles that have nonempty offspring after $Rt$
steps has Poisson distribution with parameter $(x/t^2)$.
Using this observation we can calculate the Laplace transform of the 
rescaled branching process $\xi$
\begin{eqnarray*}
\lim_{R\to\infty} \int_0^\infty e^{-\lambda y}
\P\B\{ \f{2\xi_{Rt}}{R^2} \in dy \B| \f{2\xi_0}{R^2} \B\}
&=& 
\sum_{j=0}^\infty e^{-x/t^2} \f{(x/t^2)^j}{j!} 
\B( 1 - \f{\lambda t^2}{(1+t \sqrt{\lambda})^2}\B)^j
\end{eqnarray*}
which gives
\begin{equation}\label{eq.utl}
u_t(\lambda) = \f { \lambda }{ (1 + t \sqrt\lambda)^2 }. 
\end{equation}
Thus the only possible limit for the process \eqref{xitr} is a 
continuous-state branching process $\zeta_t$
characterized by \eqref{utl}, 
and indeed, as shown by Grimvall in \cite{Gr},
such convergence holds in $D[0,1]$.

Observing that $u_t(\lambda)$ verifies 
\[
u_t(\lambda) + \int_0^t \psi(u_t(\lambda)) ds = \lambda
\]
for $\psi(u) = 2u^{3/2}$, 
we proceed with the proof of 
Theorem~\ref{T2}.

We can now state the main result of this section
\begin{theorem}\label{T3}
Let $\hat\zeta_t$ be a continuous-state branching process 
with branching mechanism $\psi(u) = 2u^{3/2}$,
started from the initial distribution $\Gamma(3/2)$ at $t=0$
and conditioned to extintion at time $t=1$.
Then the following convergence holds
\begin{equation}\label{eq.gatr2}
\f{2 |\gamma_{tR}|}{R^2} \to \hat\xi_{1-t} 
\end{equation}
as $R\to\infty$ in $D[0,1]$.
\end{theorem}

The proof consists of two parts. 
First we show that the finite-dimensional
joint distributions of the left hand side of \eqref{gatr2} 
converge to those of $\hat\zeta_{1-t}$,
then we verify the certain tightness conditions that imply convergence in
$D[0,1]$.

\subsection{Convergence of finite-dimensional distributions}
Partition $[0,1]$ into $k$ subintervals $t_1+\ldots+t_k=1$ 
and write $\{\zeta: x\tot{t} x'\}$ for the 
event "process $\zeta$ started from state $x$ and reached state $x'$ in time $t$".
Then
the multidimensional Laplace transform
of $\hat\zeta_{1-t}$ is given by
\begin{equation}\label{eq.ixik}
\int\limits_{x_1,\ldots,x_k>0}
  \f{2\sqrt{x_k} e^{-x_k}}{\sqrt{\pi}}\,
  \f{ \P\{ \zeta: dx_k \tot{t_k} dx_{k-1} \tot{t_{k-1}} \ldots
       \tot{t_2} dx_1 \tot{t_1} 0 \}}
	{ \P\{ x_k \tot{t_1+\ldots+t_k} 0 \}}
  \,\exp\B(-\sum_{j=1}^k s_j x_j\B)
\end{equation}
We wish to evaluate this expression by taking the $k$ integrals
one by one, starting from $x_1$.
First,
\[ \P\{ \zeta_t=0 | \zeta_0 = x \}
 = \lim_{\lambda\to\infty} \exp(-x u_t(\lambda)) 
 = e^{-x/t^2},
\]
and
\[ \P\{ \zeta: x_1\tot{t_1} 0 \} 
 := \lim_{dt\to 0} 
    \f{ \P\{ \zeta_{t}>0, \zeta_{t+dt}=0 | \zeta_0 = x_1 \} }{ dt }
 = \f{2x_1}{t_1^3} e^{-x_1/t_1^2}.
\]
Now fix $x_2,\ldots,x_k$ and integrate \eqref{ixik} with respect to $x_1$.
Since 
\[ \int_{0}^\infty y e^{-\lambda y} \P\{x\tot{t}dy\}
   = x u'_t(\lambda) \exp( - x  u_t(\lambda))
\]   
we have (omitting terms in \eqref{ixik} that do not depend on $x_1$)
\[ \int_0^\infty 
   \P\{ x_2 \tot{t_2} dx_1 \} e^{-s_1 x_1} x_1 e^{-x_1/t_1^2}
   = \f{x_2}{(1+t_2\sqrt{\lambda_1})^3} e^{ -x_2 u_{t_2}(\lambda_1) },
\]
where $\lambda_1 = s_1 + t_1^{-2}$.
Put $\lambda_j = s_j +  u_{t_j}(\lambda_{j-1})$ for $j=2,\ldots,k$.
The next $(k-2)$ integrals 
are calculated similarly:
\[ \int_0^\infty
   \P\{ x_{j+1} \tot{t_{j+1}} dx_j \} x_j e^{-\lambda_j x_j} 
   = \f{x_{j+1}}{(1+t_{j+1}\sqrt{\lambda_j})^3} 
      e^{ -x_{j+1} \phi_{t_{j+1}}(\lambda_j) }.
\]
Finally we'll arrive to
\begin{eqnarray}
\eqref{ixik}
&=& \int_0^\infty 
\f{x_k^{1/2}}{\sqrt\pi} e^{-x_k \lambda_k} \,dx_k
\,\cdot\,
\prod_{j=1}^{k-1} \f1{(1+t_{j+1}\sqrt{\lambda_j})^3} 
\,\cdot\,
\f{2}{t_1^3}
\nonumber\\
&=& \B[
    \lambda_k^{1/2} \cdot
	(1 + t_{k}\sqrt{\lambda_{k-1}})
	\cdot\ldots\cdot
	(1 + t_{2}\sqrt{\lambda_{1}}) \cdot
    t_1 
	\B]^{-3}
	\label{eq.lapl1}
\end{eqnarray}
where
\begin{equation}\label{eq.lrec}
 \lambda_1 = s_1 + t_1^{-2}, \qquad 
   \lambda_j = s_j +  u_{t_j}(\lambda_{j-1})
             = (\lambda_{j-1}^{-1/2}+t_j)^{-2}.
\end{equation}

Now let us calculate the limit of analogous multidimensional Laplace 
transform for the process $2 |\gamma_{Rt}|/R^2$.

Write $\{ \xi : n {\stackrel{r}\to} m \}$ for the event
"process $\xi$ started from state $n$ and arrived to 
state $m$ after $r$ steps", and assume, as eralier $t_1+\ldots+t_k=1$.

We have
\[ \P\{ \xi: n \tot{r} m \} = \f{[t^n] F(t)}{[t^m] F(t)} [t^m] \phi_r^n(t), \]
%
%
\begin{eqnarray} 
\sum_{n_1,\ldots,n_k=1}^\infty
\lefteqn{
\P\{ \xi: 
n_k \tot{r_k} n_{k-1} \tot{r_{k-1}}
\ldots \to n_2 \tot{r_2} n_1 \tot{r_1} 1
\} 
\,\,e^{-( a_1 n_1 + \ldots +  a_k n_k)}
}&&\nonumber \\
&=& 
\sum_{n_1,\ldots,n_k=1}^\infty
[t^{n_k}] F(t) \cdot 
e^{- a_k n_k} [t^{n_{k-1}}] \phi_{r_k}^{n_{k}}(t) \cdot
\ldots \nonumber \\
&&{} 
\phantom{\sum_{n_1,\ldots,n_k}^\infty}
\ldots \cdot
e^{- a_2 n_2} [t^{n_{1}}] \phi_{r_2}^{n_{2}}(t) \cdot
e^{- a_1 n_1} [t] \phi_{r_1}^{n_{1}}(t)  \cdot
\f{1}{[t]F(t)} \nonumber \\
& = & \f32\,\, 
[t] F( e^{- a_k} \phi_{r_k} 
 ( e^{- a_{k-1}} \phi_{r_{k-1}} 
 ( \ldots
 ( e^{- a_1} \phi_{r_1}(t) ) \ldots ))).
\label{eq.fold}
\end{eqnarray}
Substition $r_i = t_i R$, $ a_i = 2 s_i/R^2$ in \eqref{fold} 
gives multidimensional Laplace transform for the rescaled process,
and we want to show that as $R\to\infty$
\begin{eqnarray}
\f32
\f{\d}{\d t} \lefteqn{
\B[
F( e^{-\f{2s_k}{R^2}} \phi_{t_k R} 
 ( e^{-\f{2s_{k-1}}{R^2}} \phi_{t_{k-1} R} 
 ( \ldots 
 ( e^{-\f{2s_1}{R^2}} \phi_{t_1 R}(t) ) \ldots ))) 
 \B] \B|_{t=0}
 }
\nonumber\\
&=& 
\B[
    \lambda_k^{1/2} \cdot
	(1 + t_{k}\sqrt{\lambda_{k-1}})
	\cdot\ldots\cdot
	(1 + t_{2}\sqrt{\lambda_{1}}) \cdot
    t_1 
	\B]^{-3} + O(1)
\label{eq.flap}
\end{eqnarray}
with $\lambda_1,\ldots,\lambda_k$ defined by \eqref{lrec}.

We'll start with expanding the innermost instance of $\phi$ in \eqref{flap},
$\phi_{t_1R}$, into series in~$R$ at $t=0$.
\[ \phi_{t_1 R}(0) = 1 - \f{2}{t_1 ^2 R^2} + O(\f1{R^3}), \qquad
   e^{-\f{2s_1}{R^2}} = 1 - \f{2s_1}{R^2} + O(\f1{R^3}). \]
Put $x_1 = s_1 + t_1^{-2}$.
Now the argument to $\phi_{t_2R}$ in \eqref{flap}
is $1 - 2x_1/R^2 + O(1/R^3)$.
Further expansion gives
\[ \phi_{t_j R}\B( 1 - \f{2x_{j-1}}{R^2} + O(\f1{R^3}) \B)
   = 1 - \f{2}{(x_{j-1}^{-1/2} + t)^2 R^2} + O(\f1{R^3}), \]
so putting $x_j = s_j + (x_{j-1}^{-1/2} + t_j)^{-2}$
we can continue expanding $\phi_{t_jR}$ one by one.

Finally
\[ \phi'_{t_1 R}(0) = \f83\,\f1{t_1^3 R^3} + O(\f1{R^4}), \]
\[ F'\B(1 - \f{2 x_k}{R^2} + O(\f1{R^3})\B)
   = \f14\,\f{R^3}{x_k^{3/2}} + O(R^2), 
\]   
\[ \phi'_{t_j R}\B( 1 - \f{2x_{j-1}}{R^2} + O(\f1{R^3}) \B)
   = \f1{(1 + t_j x_{j-1}^{1/2})^3} + O(\f1R),
\]
and
\begin{eqnarray*}
\f32
\f{\d}{\d t} 
\lefteqn{ 
\B[
F( e^{-\f{2s_k}{R^2}} \phi_{t_k R} 
 ( e^{-\f{2s_{k-1}}{R^2}} \phi_{t_{k-1} R} 
 ( \ldots 
 ( e^{-\f{2s_1}{R^2}} \phi_{t_1 R}(t) ) \ldots ))) 
 \B] \B|_{t=0}
 }
\nonumber\\
&=&
\f32\, 
 F'( e^{-\f{2s_k}{R^2}} \phi_{t_k R} 
   ( e^{-\f{2s_{k-1}}{R^2}} \phi_{t_{k-1} R} 
   ( \ldots 
   ( e^{-\f{2s_1}{R^2}} \phi_{t_1 R}(0) ) \ldots )))  \nonumber\\
&& \phantom{ f\,( e^{-\f{2s_k}{R^2}} } \times \phi'_{t_k R} 
   ( e^{-\f{2s_{k-1}}{R^2}} \phi_{t_{k-1} R} 
   ( \ldots 
   ( e^{-\f{2s_1}{R^2}} \phi_{t_1 R}(0) ) \ldots ))   \nonumber\\
&& \phantom{f\,( e^{-\f{2s_k}{R^2}} \phi_{t_k R} 
   ( e^{-\f{2s_{k-1}}{R^2}} } \times \phi'_{t_{k-1} R} 
   ( \ldots 
   ( e^{-\f{2s_1}{R^2}} \phi_{t_1 R}(0) ) \ldots )    \nonumber\\
&& \phantom{f\,( e^{-\f{2s_k}{R^2}} \phi_{t_k R} 
   ( e^{-\f{2s_{k-1}}{R^2}} \phi'_{t_{k-1} R} 
   ( \ldots ( e^{-\f{2s_1}{R^2}} \times} 
   \ldots
   \nonumber\\
&& \phantom{f\,( e^{-\f{2s_k}{R^2}} \phi_{t_k R} 
   ( e^{-\f{2s_{k-1}}{R^2}}  \phi'_{t_{k-1} R} 
   ( \ldots ( e^{-\f{2s_1}{R^2}} } 
 \times \phi'_{t_1 R}(0) 
 \vphantom{\int_0}
 \nonumber\\
&=&
   \B[ x_k^{1/2} (1 + t_k x_{k-1}^{1/2}) 
                 (1 + t_{k-1} x_{k-2}^{1/2}) \cdots
                 (1 + t_2 x_1^{1/2}) t_1 \B]^{-3}
				 + O(1)
\end{eqnarray*}
where
\[ x_1 = s_1 + t_1^{-2},\qquad
   x_j = s_j + (x_{j-1}^{-1/2} + t_j)^{-2}. \]
which is exactly the same expression as \eqref{lapl1}.


\subsection{Tightness}
The main tool in proving tightness will be the Theorem $2.2'$ in \cite{Gr}.
\begin{theorem}[Grimvall]\label{TG}
Let $\{X_k^{(n)}\}_{n=1}^\infty$ be a sequence of real-valued Markov chains
and let
$\nu_a^{(n)}$ denote a measure defined by 
\[ \nu_a^{(n)} = \P\{ X^{(n)}_{k+1} \in E | X^{(n)}_k = a \} \]
for all $a\in \mathbb{R}$ and all Borel sets E.
Let also 
\[ Y^{(n)}(t) = X^{(n)}_{\lceil nt \rceil}. \]
Then the sequence $\{Y^{(n)}\}_n$ is tight in $D[0,1]$, if
\begin{list}{}{}
\item{\rm(i)} $\P\{ \sup_{0\le t\le 1} |Y^{(n)}(t)| > \lambda \} \to 0$
 as $\lambda\to\infty$ uniformly in $n\in\mathbb{N}$,
\item{\rm(ii)} $\{(\nu_a^{(n)})^{*n}\}_{a\in C, n\in \mathbb{N}}$
is tight for every compact subset $C$ of the real line.
\end{list}
\end{theorem}
\vskip.5cm
\noindent
We wish to apply this theorem with
\[ X_k^{(R)} = \f{2 \xi_k}{R^2}, \qquad
   Y^{(R)}(t) = \f{2 \xi_{\lceil tR \rceil}}{R^2}. \]
The measure $\nu_a^{(R)}$ then corresponds to the conditional distribution
\begin{equation}\label{eq.finr}
\B( \f{\xi_1 - m}{R^2/2} \B| \xi_0 = m \B), 
   \quad m=aR^2/2. 
\end{equation}
Using the representation
\[ \E( \xi_1 | \xi_0 = m) = \f{ [t^m] F'(\phi(t)) }{ [t^m] F(t) }
\quad
   \E( \xi_1 (\xi_1 - 1) | \xi_0 = m) = \f{ [t^m] F''(\phi(t)) }{ [t^m] F(t) }
\]
and expanding the functions $F$, $F'(\phi)$, $F''(\phi)$ in series near $t=1$
(thanks to the explicit formulas for $\phi$ and $F$ this calculation becomes 
trivial)
we find the asymptotics
\begin{equation}\label{eq.exim1}
 \E(\xi_1 - m |\xi_0=m) = \f34 \sqrt{2\pi} m^{1/2} + O(1) 
\end{equation}
\begin{equation}\label{eq.exim2}
 \E( (\xi_1 - m)^2 |\xi_0=m) = \f38 \sqrt{2\pi} m^{3/2} + O(m^{1/2}). 
\end{equation}

It follows from \eqref{exim1} that $\xi_R$ is a submartingale
and by Kolmogorov-Doob inequality
\[ \P\{ \sup_{0\le k\le R} \xi_k \ge \lambda R^2 \}
   \le \f{  \E \xi_R }{\lambda R^2}   
   \to 0 \quad \mbox{as $\lambda\to\infty$},
\]
thus the condition (i) of Theorem~\ref{TG} is satisfied.

Next 
note that the $n$-fold convolution of 
$\nu_a^{(n)}$ corresponds to the sum of $n$ independent copies of \eqref{finr},
and that by \eqref{exim1}, \eqref{exim2} this sum has both mean and
variance of order $O(1)$ as $R\to\infty$.
The condition (ii) of Theorem~\ref{TG} is then satisfied by e.g. Chebyshev 
inequality.

Thus the sequence of processes \eqref{gatr} is tight in $D[0,1]$,
and by Prohorov theorem the convergence \eqref{gatr2} follows.
This finishes the proof of Theorem~\ref{T3}.

%


\section{Enumeration}\label{sec.p}
The formula~\eqref{qx} is obtained from a more general formula for the number of
bicubic (bipartite, trivalent) planar maps due to Tutte \cite{T1}.
No doubt, \eqref{Uxy} could also be derived from the same source, 
but we shall give a slightly more straightforward proof.

Consider first the class $\calQ'$ of quadrangulations with simple boundary 
with no double edges. Every such quadrangulation has at least two faces.
Take the vertex opposite to the root vertex in the rooted face
and cut the quadrangulation along every edge, incident to this vertex.
Forget the rooted face.
This operation produces one or more components, each being either a single face,
either again a quadrangulation from $\calQ'$, and the boundary of each component
consists of two segments -- one is a part of original quadrangulation's
boundary, another consists of previously internal faces (see \figref{req}).

\putfigure{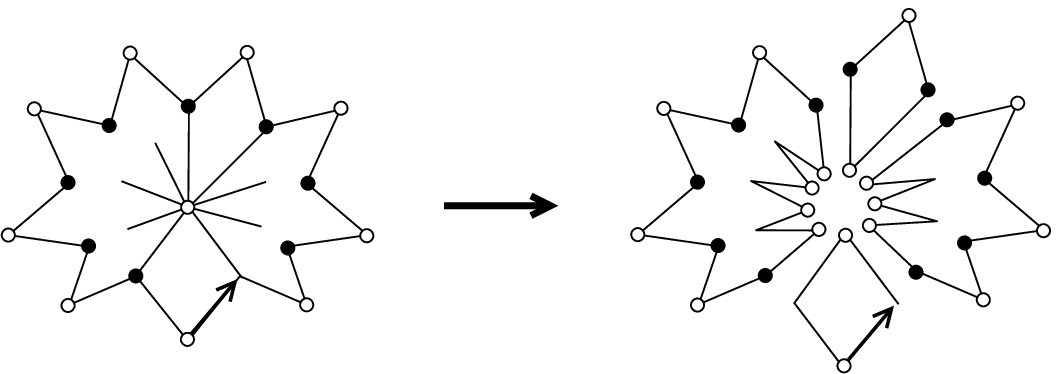}{}

This is a bijection -- given an ordered collection elements, each being either
a single face, either a quadrangulations with boundary separated in two 
segments of length at least $1$, 
the original quadrangulation can be reconstructed.
Thus
\[ u(x,y) = xy \cdot \sum_{n=1}^\infty \B( xy + w(x,y,1) \B)^n, \]
where $u(x,y)$ is the generating function of quadrangulations in $\calQ'$ with 
$x$ counting faces and $y$ counting faces at the boundary, and
\[ w(x,y,z) = \f{ y\, u(x,z) - z\, u(x,y) }{ z - y } \]
is the generating function of quadrangulations with segmented boundary, $y$ and
$z$ counting the length (in faces) of each segment.

Expressing $w(x,y,1)$ from the first equation 
and substituting into the second one with $z=1$, we get
\[ 1 - xy - \f{xy}{u(x,y) + xy} = \f{ y u(x,1) - u(x,y) }{ 1-y }.  \]
Note that this equation is quadratic in $u(x,y)$.

It's not hard to see that $u(x,1)$ is the generating function of 
spere quadrangulations without double edges and {\em with at least $2$ faces}.
There are two possible quadrangulations with a single face, 
and one more special quadrangulation with zero faces,
so the complete generating function of quadrangulations without 
double edges is
\[ u(x,1) + 2x + 1. \]


Now to pass from quadrangulations without double edges to the general 
class of quadrangulations $\calQ$ we attach at each internal edge of
quadrangulation from $\calQ'$ a general sphere quadrangulation. 
More precisely, we cut this edge and
identify two sides of obtained hole with two sides of analogous hole, obtained
by cutting the root edge of quadrangulation being attached.
This is the \defined{extension} procedure, best explained in \cite{T1}, 
section 7.

Etension is equivalent to the substitution $x\to q^2x$, $y\to q^{-1}y$.
Under this substitution we get
\[ u( q^2 x, 1) + 2x q^2 + 1 = q, \]
\[ u( q^2 x, q^{-1}y) + xy\,q = U(x,y). \]
The term $xyq$ in the second equation corresponds to the quadrangulations
with simple boundary of length $1$.

Combining two last equations with quadratic equation on $u(x,y)$ we get
\eqref{Uxy}.

\section{Discussion}
Random planar maps are  considered a natural model of space with
fluctuating geometry in 2-dimensional quantum gravity.
Ambjorn and Watabiki \cite{AW} suggested that the internal Hausdorf dimension
such random space is $4$, and this relation doesn't depend on the  choice of
triangulations, quadrangulations, or some other reasonable distribution of
polygons as the underlying model.

Theorem~\ref{T1} for triangulations was proved by Angel and Schramm
\cite{AS}. They also provided estimates for the growth of hull boundary,
which was shown to be quadratic in radius with polylog corrections \cite{A},
but no exact limit was found.

A similar theorem was proved by Chassaing and Durhuus~\cite{CD},
who showed the local weak convergence of well-labeled trees, which are known
to be bijective to quadrangulations~\cite{ST}. 
This bijection, however, is continuous only in one direction: from
quadrangulations to trees (with respect topology of $\calQ$ and the natural
local topology on the space of trees),
so this result is not equivalent to Theorem~\ref{T1}.

Related model of random triangulations with free boundary was considered by
Malyshev and Krikun \cite{MK}. It was shown that at the critical boundary
parameter value the boundary of a random quadrangulation with $N$ faces has
about $\sqrt{N}$ edges and the ratio converges in distribution to the square of
$\Gamma(3/4)$ law. Nothing is known about the diameter of the triangulation, 
but this is natural to suggest that the diameter has order $N^{1/4}$.

The skeleton construction used in this paper was proposed in~\cite{K},
where it was applied to random triangulations. The branching process obtained
for triangulations differs from $\xi$, 
but it has non-extinction probabilities of the same order $1/R^2$
and it's generating function has the main singularity of the same order $3/2$.


There should also exist certain natural mapping of branching process structure
into the brownian map \cite{MM}.

Finally, we want to note the following statement in \cite{AW}:
{\em "A boundary of $l$ links will have the discrete length $l$ in lattice
units, but if we view the boundary from the interior of the surface its true
linear extension $r$ will only be $\sqrt{l}$, since the boundary can be viewed
as a random walk from the interior".}

\vskip 3cm
\noindent
Maxim Krikun \\
Institut Elie Cartan, \\
Universite Henri Poincare,\\
Nancy, France\\
krikun@iecn.u-nancy.fr \\

\end{document}

%% file: pre.tex
\def\E{\mathop{\hbox{\sf E}}\nolimits}
\def\P{\mathop{\hbox{\sf P}}\nolimits}

\def\phi{\varphi}

\def\f{\frac}
\def\d{\partial}

\def\B{\Big}

\newtheorem{lemma}{Lemma}[section]
\newtheorem{theorem}{Theorem}

\newtheorem{conjecture}{Conjecture}

\usepackage{graphics}

\def\figref#1{fig. \ref{fig.#1}}
\def\eqref#1{(\ref{eq.#1})}

\def\putfigure#1#2{
	\begin{figure}[ht]
	\centering
	\includegraphics{#1.eps}
	\caption{#2}
	\label{fig.#1}
	\end{figure}
}

\def\defined#1{{\em #1}}

\def\optional#1{}

\def\intl{\int\limits}

\def\proof{\par\noindent{\bf Proof.\ }}
\def\proofof#1{\par\noindent{\bf Proof of #1.\ }}
\def\eop{\vskip 3mm }